\documentclass[12pt,a4paper]{article}

\usepackage{amssymb,amsmath,amsfonts}
\newcommand{\bb}{\begin{equation}}
\newcommand{\ee}{\end{equation}}

 \newtheorem{te}{Theorem}[section]

 \newtheorem{de}{Definition}[section]
\newcommand{\QED} {\hfill$\square$}

   \font\sst=cmtt8
  \font\ssi=cmti8 
  \font\sst=cmtt8 
 \font\ssi=cmti8

\def\PXX#1#2#3{ P^{#1}_{#2#3}}

\def\D#1#2{\delta^{#1}_{#2}\,}
\def\KX#1#2#3{\xi^{#1}_{#2#3}\,}

\def\GN#1#2{\overline g_{#1#2}\,}

\def\AA#1#2#3#4#5{\underset {#1} A^{#2}_{\,\,#3#4#5}}

\def \AK#1{a^i_{j\underset {#1} |m}}
\def \AZ{a^i_{j,m}}
\def \AA#1#2{a^{#1}_{#2}}

\def\PXX#1#2#3{ P^{#1}_{#2#3}}

\def\D#1#2{\delta^{#1}_{#2}\,}
\def\KX#1#2#3{\xi^{#1}_{#2#3}\,}

\def\GN#1#2{\overline g_{#1#2}\,}

\def\KK#1#2#3{\xi^{#1}_{#2#3}}

\def\GA#1#2#3{ \Gamma^{#1}_{\underset \lor {#2#3}}}

\def\GG#1#2#3{ \Gamma {}^{#1}_{#2#3}}

\def\PXX#1#2#3{ P^{#1}_{#2#3}}

\def\D#1#2{\delta^{#1}_{#2}\,}
\def\KX#1#2#3{\xi^{#1}_{#2#3}\,}

\def\GN#1#2{\overline g_{#1#2}\,}

\def\AA#1#2#3#4#5{\underset {#1} A^{#2}_{\,\,#3#4#5}}

\def\PXX#1#2#3{ P^{#1}_{#2#3}(u)}

\def\D#1#2{\delta^{#1}_{#2}\,}
\def\KX#1#2#3{\xi^{#1}_{#2#3}\,}

\def\GN#1#2{\overline g_{#1#2}{}}

\def\PXX#1#2{\underset{#1}\psi{}_{#2}}

\def\KK#1#2#3{\xi^#1_{#2#3}}

\def\enddemo{\qed \endtrivlist}
\expandafter\let\csname enddemo*\endcsname=\enddemo
\def\qedsymbol{\ifmmode\bgroup\else$\bgroup\aftergroup$\fi
  \vcenter{\hrule\hbox{\vrule height.6em\kern.6em\vrule}\hrule}\egroup}
\def\qed{\ifmmode\else\unskip\nobreak\fi\quad\qedsymbol}

\def\GA#1#2#3{ \Gamma^#1_{\underset \lor {#2#3}}}

\def\GG#1#2#3{ \Gamma ^#1_{#2#3}}

\def\ff#1#2{F^{#1}_{#2}}

\def\gg#1#2{g_{\underline{#1#2}}{}}
\def\ggg#1#2{g^{\underline{#1#2}}{}}

\def\PXX#1#2#3{ P^{#1}_{#2#3}}

\def\D#1#2{\delta^{#1}_{#2}\,}
\def\KX#1#2#3{\xi^{#1}_{#2#3}\,}

\def\GN#1#2{\overline g_{#1#2}\,}

\def\AA#1#2#3#4#5{\underset {#1} A^{#2}_{\,\,#3#4#5}}

\def\PXX#1#2#3{ P^{#1}_{#2#3}}

\def\GN#1#2{\overline g_{#1#2}\,}

\def\KK#1#2#3{\xi^#1_{#2#3}}
\def\GA#1#2#3{ \Gamma^{*#1}_{\underset \lor {#2#3}}}

\def\GG#1#2#3{ \Gamma ^#1_{#2#3}}

\def\ff#1#2{F^{#1}_{#2}}

\def\gg#1#2{g_{#1#2}}
\def\ggs#1#2{g_{\underline{#1#2}}}
\def\gggs#1#2{g^{\underline{#1#2}}}
\def\ggg#1#2{g^{{#1#2}}}

\def \GN#1#2#3{\overline\Gamma {}^{#1}_{#2#3}}
\def \AK#1{a{}^i_{j\underset #1|m}}
\def \AZ{a^i_{j,m}}
\def \AA#1#2{a{}^{#1}_{#2}}

\def \AK#1{a^i_{j\underset #1 |m}}
\def \AZ{a^i_{j,m}}
\def \AA#1#2{a^{#1}_{#2}}

\def\AA#1#2#3#4#5{\underset {#1} A{}^{*#2}_{\,#3#4#5}}

\def \AK#1{a^i_{j\underset #1 |m}}
\def \AZ{a^i_{j,m}}
\def \AA#1#2{a^{#1}_{#2}}

\def \GG#1#2#3{\Gamma ^{#1}_{#2#3}}

\def\PXX#1#2#3{ P^{#1}_{#2#3}}
\def\D#1#2{\delta^{#1}_{#2}\,}
\def\KX#1#2#3{\xi^{#1}_{#2#3}\,}

\def\KK#1#2#3{\xi^#1_{#2#3}}
\def\GA#1#2#3{ \Gamma^#1_{\underset \lor {#2#3}}}

\def\ff#1#2{F^{#1}_{#2}}
\def\gg#1#2{g_{#1#2}}
\def\ggg#1#2{g^{#1#2}}
\def\ggs#1#2{g_{\underline{#1#2}}}
\def\gggs#1#2{g^{\underline{#1#2}}}

\def\ff#1#2{F^{#1}_{#2}}

\def\GA#1#2#3{ \Gamma{}^{#1}_{\underset \lor {#2#3}}}

\def\gs#1#2{g{}_{\underline {#1#2}}\,}
\def\ggs#1#2{g^{\underline {#1#2}}\,}
\def\gn#1#2{g{}_{\underset \lor  {#1#2}}\,}

  \font\sst=cmtt8 
  \font\ssi=cmti8

\date{}

\title{\bf Geodesic mapping onto K\"ahlerian space of the first kind\footnote{The authors gratefully acknowledge
support from the research projects 174012 of the Serbian Ministry of
Science and  FAST-S-12-25/1660 of the Brno University of Technology.}}

\author{\frenchspacing Milan Zlatanovi\'c,$^1$\;  Irena Hinterleitner$^2$ and Marija Najdanovi\' c$^3$ \\
{\ssi $^{1,3}$ Faculty of Sciences and Mathematics, University of Ni\v{s}, Serbia}\\
{\ssi $^2$ Department of Mathematics, Faculty of Civil Engineering, Brno University of Technology, Czech Republic} \\
{\ssi E-mails:}  {\sst $^1$ zlatmilan@pmf.ni.ac.rs},\ {\sst $^2$ Hinterleitner.Irena@seznam.cz}, \ {\sst $^3$ marijamath@yahoo.com}\\
}

\begin{document}
\maketitle

\begin{abstract}
In the present paper a generalized K\"
ahlerian space $\mathbb{G}\underset 1 {\mathbb{K}}{}_N$ of the first kind is considered,
as a generalized Riemannian space $\mathbb{GR}_N$ with almost complex structure $F^h_i$, that is
covariantly constant with respect to the first kind of covariant
derivative.

 Using the non-symmetric metric tensor
we find necessary and sufficient conditions for a geodesic mapping $f:\mathbb{GR}_N\to \mathbb{G}\underset 1 {\mathbb{\overline{K}}}{}_N$
 with respect to the four kinds of covariant derivatives.

\begin{description}
\item[] {\it AMS Subj. Class.:} 53B05, 53B35.
\item[] {\it Key words:} geodesic mapping,
generalized K\"ahlerian space, equitorsion geodesic mapping.
\end{description}
\end{abstract}

\section{Introduction}

An $N$-dimensional Riemannian space with metric tensor $g_{ij}$
is a \emph{K\" ahlerian space} $\mathbb{K}_N$ if there exists an
almost complex structure  $\ff ij$, such that
\begin{eqnarray}\aligned
\ff hp\ff pi&=-\D hi,\\
\gg pq\ff pi\ff qj=\gg ij&,\quad \ggg ij=\ggg pq\ff ip\ff jq,\\
F^h_{i;j}&=0,\quad \endaligned\end{eqnarray}
 where $(;)$  denotes the covariant derivative with respect to the 
metric tensor  $\gg ij$.
\medskip

 A \emph{generalized Riemannian space} $\mathbb{GR}_N$ in the sense of
Eisenhart's definition \cite {bib-1} is a differentiable
$N$-dimensional manifold, equipped with a non-symmetric basic tensor
$g_{ij}$ (i.e. $g_{ij}\ne g_{ji}$).  Based on non-symmetry, the
symmetric part and antisymmetric part of $g_{ij}$ can be defined
\begin{equation}\label{1.2}
\gs ij=\frac{1}{2}(g_{ij}+g_{ji})=\frac{1}{2}g_{(ij)},\quad \gn
ij=\frac{1}{2}(g_{ij}-g_{ji})=\frac{1}{2}g_{[ij]}.
\end{equation}
 The lowering and the raising of indices is defined by the
tensors $\gs ij$ and $\ggs ij$ respectively, where  $\ggs ij$ is
defined by the equation
\begin{equation}\label{1.4}
{\gs ij\ggs jk=\delta_i^k}
\end{equation}
and $\delta_i^k$ is the Kronecker symbol. From (\ref{1.4}) we have that
the matrix ($\ggs ij$) is inverse to $(\gs ij)$, wherefrom  it is
necessary to be $ g=det(\gs ij)\neq0. $ Connection coefficients of
this space are generalized Cristoffel symbols of the second kind,
where
\begin{equation}\label{1.5}
{\Gamma_{i.jk}=\frac{1}{2}(g_{ji,k}-g_{jk,i}+g_{ik,j}),}\qquad
g_{ij,k}=\frac{\partial g_{ij}}{\partial x^k}.
\end{equation}
\begin{equation}\label{1.6}
{\Gamma_{jk}^i=\ggs ip\Gamma_{p.jk}.}
\end{equation}
Generally $\Gamma ^i_{jk}\ne \Gamma ^i_{kj}$. Therefore,
one can define the symmetric part and anti-symmetric part of
$\Gamma^i_{jk}$, respectively
\begin{equation}
\Gamma^i_{\underline{jk}}=\frac
12(\Gamma^i_{jk}+\Gamma^i_{kj})=\frac 12 \Gamma^i_{(jk)},\quad
\Gamma^i_{\underset\lor {jk}}=\frac
12(\Gamma^i_{jk}-\Gamma^i_{kj})=\frac 12 \Gamma^i_{[jk]}
\end{equation}%
The quantity $\Gamma^i_{\underset\lor{jk}}$ is the {\it torsion tensor}
of the spaces $\mathbb{GR}_N$.
\medskip

The use of a non-symmetric metric tensor and a non-symmetric connection
became especially topical after the appearance of the papers of A.
Einstein \cite{a1}-\cite{a4} related to the creation of the Unified
Field Theory (UFT). More recently the idea of a non-symmetric metric tensor appears in Moffat's non-symmetric gravitational theory \cite{mof}. 
We remark that in the UFT the symmetric part $
g_{\underline{ij}}$ of the metric tensor $g_{ij}$ is related to 
gravitation, and  antisymmetric one $g_{\underset\lor{ij}}$ to
the electro\-magnetism. In Moffat's theory the antysymmetric part represents a Proca field (massive Maxeell field) which is part of the gravitational interaction, contributing to the rotation of galaxies \cite{mof}.

Based on non-symmetry of the connection in  a generalized Riemannian space one
can define four kinds of covariant derivatives.
For example, for a tensor $a^i_j$ in $\mathbb{GR}_N$ we have
\begin{eqnarray}\aligned
 &\AK 1=\AZ +\GG {i}{p} m\AA {p} {j}-\GG p jm\AA ip ,
\\
&\AK 2=\AZ +\GG imp \AA p j-\GG p mj\AA ip ,\\
 &\AK 3=\AZ +\GG ip m\AA p
j-\GG p mj\AA ip ,
 \\
 &\AK 4=\AZ +\GG imp \AA p j-\GG p jm\AA ip .\endaligned
\end{eqnarray}

By applying four kinds of covariant derivatives of a tensor, it is
possible to construct several Ricci type identities. In these
identities 12 curvature tensors appear and 15 quantities which are
not tensors, named "curvature pseudotensors" by Min\v ci\' c
\cite{bib-2,bib-3}.
 In the
case of the space $\mathbb{GR}_N$ we have five independent curvature tensors.


\subsection{Generalized K\"ahlerian space of the first kind}

 K\" ahlerian spaces and their mappings were investigated by many
authors, for example K. Yano, M. Prvanovi\'c,
 J. Mike\v s,
Doma\v sev, N. Pu\v si\'c, T. Otsuki and Y. Tasiro, S. S. Pujar,
 N.S. Sinjukov and many other
authors.

\begin{de} A generalized $N$-dimensional Riemannian
space with non-symmetric metric tensor $\gg ij$, is a {\bf
generalized K\" ahlerian space of the first kind} $\mathbb{G}\underset 1
{\mathbb{K}}{}_N$  if there exists an almost complex structure
$\ff ij$, so that
\begin{equation}\ff hp\ff pi=-\D hi,\label{hpv2.1}\end{equation}
\begin{equation}\gs pq\ff pi\ff qj=\gs ij,\quad \gggs ij=\gggs pq\ff ip\ff
jq,\label{hpv2.2}\end{equation}
\begin{equation}F^h_{i\underset1 |j}=0,\label{hpv2.3}\end{equation}
\begin{equation}F^h_{i;j}=0,\label{hp3}\end{equation}
where $\underset 1 |$  denotes the covariant derivative of
the first kind with respect to the  connection $\Gamma^i_{jk}$ $(\Gamma^i_{jk}\ne\Gamma^i_{kj})$
 and $(;)$ denotes the covariant derivative with respect to the  symmetric part of the metric tensor $\Gamma^i_{\underline{jk}}$. \end{de}

From
(\ref{hpv2.2}), using  (\ref{hpv2.1}), we get $F_{ij}=-F_{ji}, \;
F^{ij}=-F^{ji},$ where we denote $ F_{ji}=\ff pj\ggs pi,\;
F^{ji}=\ff jp\gggs pi$.

\medskip
From this follows the theorem:

\begin{te} For the almost complex structure $\ff ij$ of  $\mathbb{G}\underset 1
{\mathbb{K}}{}_N$ the  relations
\begin{equation}
F^h_{i\underset 2 |j}=0,\quad
F^h_{i\underset 3 |j}=2\ff hp\GA pij,\quad F^h_{i\underset 4
|j}=2\ff pi\GA hjp\label{hpv2.4}
\end{equation}
are valid, where $\GA hij$ is a torsion tensor.\end{te}

\setcounter{equation}{0}
\section{Geodesic mapping}

In this part we consider a geodesic mapping $f:\mathbb{GR}_N \to
\mathbb{G}\overline {\underset 1{\mathbb{K}}}{}_N.$

\begin{de}A diffeomorphism $f:\mathbb{GR}_N \to
\mathbb{G}\overline {\underset 1{\mathbb{K}}}{}_N$ is {\bf geodesic}, if geodesics of the space $\mathbb{GR}_N$ are mapped to geodesics of
the space $\mathbb{G}\overline {\underset1{\mathbb{K}}}{}_N$.
\end{de}

 In the corresponding
points $M$ and $\overline M$ we can put
\begin{equation}
\overline{\Gamma}{} ^i_{jk}=\Gamma ^i_{jk}+\PXX ijk,\quad
(i,j,k=1,...,N),\label{0.2}
\end{equation}
where $\PXX ijk$ is {\it the deformation tensor\/} of the connection
$\Gamma$ of $\mathbb{GR}_N$ according to the mapping $f:\mathbb{GR}_N \to
\mathbb{G}\overline {\underset 1{\mathbb{K}}}{}_N$.
\medskip

\begin{te} {\bf \cite{bib-65}} A necessary and sufficient condition that the mapping $f:\mathbb{GR}_N \to
\mathbb{G}\overline {\underset 1{\mathbb{K}}}{}_N$ be
geodesic is that the deformation tensor $P^i_{jk}$
from $(\ref{0.2})$ has the form
\begin{equation}
\PXX ijk=\D i{j}\psi _{k}+\D i{k}\psi _{j}+\KK ijk,\label{0.3}
\end{equation}
where
\begin{equation}\label{0.5}
\psi _i=\frac 1{N+1}(\GN \alpha i\alpha -\GG \alpha i\alpha),\quad
\xi ^i_{jk}=P^i_{\underset \lor  {jk}}=\frac 12 (P^i_{jk}-P^i_{kj}).
\end{equation}
\end{te}
We remark that in $\mathbb{G}{\underset1{\mathbb{K}}}{}_N$ are valid the following equations:
\begin{equation}
\GA \alpha i\alpha =0,\quad \KX \alpha i\alpha =0,\quad
F^\alpha_\alpha =0.
\end{equation}
\medskip
In \cite{bibknjiga} Mike\v s at al. have proved a necessary and
sufficient condition for a geodesic mapping of a Riemannian space onto
a K\"ahlerian  space.

\begin{te} The Riemannian space
$\mathbb{R}_N$ admits a nontrivial geodesic mapping onto the 
K\"ahlerian space $ \overline{\mathbb{K}}_N$ if and only if, in the common
coordinate system $x$ with respect to the mapping,  the conditions
\begin{equation}\label{sinj}
\begin{array}{crl}
a)\;& \overline{g}_ {{{ij};k}}{}=&2\psi_k\overline{g}_{{ij}}+\psi_i\overline{g}_{{jk}}+\psi_j\overline{g}_{{ik}};\\
b)&\;\overline{F}^h_{i;k}=&\overline{F}^h_{k}\psi_i-\delta^h_k\overline{F}^\alpha_{i}\psi_\alpha;\\
\end{array}%
\end{equation}
hold, where $\psi_i\ne0$ and the tensors $\overline{g}_{{ij}}$
and $\overline{F}^h_{i}$ satisfy the following conditions:
\begin{equation}det(\overline{g}_{{ij}})\ne0,\quad
\overline{F}^h_{\alpha}\overline{g}_{{\alpha
j}}+\overline{F}^\alpha_{j}\overline{g}_{{\alpha i}}=0,
\quad
\overline{F}^h_{\alpha}\overline{F}^\alpha_{i}=-\delta^h_i.\end{equation}
Then $\overline{g}_{{ij}}$ and $\overline{F}^h_{i}$ are
the metric tensor and the structure of $ \overline{\mathbb{K}}_N$,
respectively.
 \end{te}

Our idea is to find corresponding equations with respect to the four kinds of covariant derivative.
\begin{te} The generalized Riemannian space
$\mathbb{GR}_N$ admits a nontrivial geodesic mapping onto the generalized
K\"ahlerian space $\mathbb{G}\overline{\underset1{\mathbb{K}}}{}_N$ if and only if, in
the common coordinate system $x$ with respect to the mapping,  the
conditions
\begin{equation}\label{prvav}
\begin{array}{crl}
a)\;& \overline{g}_ {{{ij}\underset1|k}}{}=&\overline{g}_
{{\underset\lor{ij}\underset1{\overline{|}}k}}{}
+2\psi_k\overline{g}_{{ij}}+\psi_i\overline{g}_{{jk}}+\psi_j\overline{g}_{{ik}}+\xi^\alpha_{ik}\overline{g}_{\alpha
j}
+\xi^\alpha_{jk}\overline{g}_{i\alpha};\\
b)&\;\overline{F}^h_{i\underset1|k}=&\overline{F}^h_{k}\psi_i
-\delta^h_k\overline{F}^\alpha_{i}\psi_\alpha
 -\xi^h_{\alpha k}\overline{F}^\alpha_{i}+\xi^\alpha_{ik}\overline{F}^h_{\alpha};
\end{array}%
\end{equation}
hold with respect to the first kind of covariant derivatives,
 where $\psi_i\ne0$ and the tensors $\overline{g}_{{ij}}$ and $\overline{F}^h_{i}$ satisfy the following conditions:
\begin{equation}\label{uslov1m}det(\overline{g}_{\underline{ij}})\ne0,\quad
\overline{F}^\alpha_{i}\overline{g}_{\underline{\alpha
j}}+\overline{F}^\alpha_{j}\overline{g}_{\underline{\alpha i}}=0,
\quad
\overline{F}^h_{\alpha}\overline{F}^\alpha_{i}=-\delta^h_i.\end{equation}
Then $\overline{g}_{{ij}}$ and $\overline{F}^h_{i}$ are
the metric tensor and the structure of
$\mathbb{G}\overline{\underset1{\mathbb{K}}}{}_N$, respectively.
 \end{te}

\begin{demo} The equation ($\ref{prvav}a$) guarantees the
existence  of geodesic mapping from a generalized Riemannian space
$\mathbb{GR}_N$  onto generalized Riemannian space
$\mathbb{G}\overline{\mathbb{R}}_N$ with metric tensor
$\overline{g}_{ij}$ with respect to the first kind of covariant
derivatives.

The formula ($\ref{prvav}b$) implies that the structure
$\overline{F}^h_{i}$ in $\mathbb{G}\overline{\mathbb{R}}_N$ is
covariantly constant with respect to the first kind of covariant
derivative. The algebraic conditions (\ref{uslov1m}) guarantee that
$\overline{g}_{ij}$ and $\overline{F}^h_{i}$ are the metric tensor
and the structure of
$\mathbb{G}\overline{\underset1{\mathbb{K}}}{}_N$, respectively.

The deformation tensor is determined by equation (\ref{0.3}), i.e.
\begin{equation}\label{def}
\overline{\Gamma}^h_{ij}-\Gamma^h_{ij}=\psi_i\delta^h_j+\psi_j\delta^h_i+\xi^h_{ij}.
\end{equation}
For the structure $\overline{F}$, we have the following equations:
\begin{equation}\label{kovar}
\overline{F}^h_{i\underset1|k}=\overline{F}^h_{i,k}+\Gamma^h_{pk}
\overline{F}^p_{i}-\Gamma^p_{ik}\overline{F}^h_{p},\quad
\overline{F}^h_{i\underset2|k}=\overline{F}^h_{i,k}+\Gamma^h_{kp}
\overline{F}^p_{i}-\Gamma^p_{ki}\overline{F}^h_{p}.
\end{equation}
Replacing $\Gamma^h_{ij}$ from (\ref{def}) in (\ref{kovar}), we get

\begin{eqnarray}\aligned
\overline{F}^h_{i\underset1|k}&=\overline{F}^h_{i,k}+(\overline{\Gamma}^h_{pk}-\psi_p\delta^h_k-\psi_k\delta^h_p-\xi^h_{pk})
 \overline{F}^p_{i}-(\overline{\Gamma}^p_{ik}-\psi_i\delta^p_k-\psi_k\delta^p_i-\xi^p_{ik})\overline{F}^h_{p}\\
 &=\overline{F}^h_{i,k}+\overline{\Gamma}^h_{pk}\overline{F}^p_{i}-\psi_p\delta^h_k\overline{F}^p_{i}
 -\psi_k\delta^h_p\overline{F}^p_{i}-\xi^h_{pk}\overline{F}^p_{i}
 -\overline{\Gamma}^p_{ik}\overline{F}^h_{p}+\psi_i\delta^p_k\overline{F}^h_{p}+\psi_k\delta^p_i\overline{F}^h_{p}+\xi^p_{ik}\overline{F}^h_{p}\\
 &=\overline{F}^h_{i\underset1{\overline{|}}k}-\psi_p\delta^h_k\overline{F}^p_{i}-\psi_k\delta^h_p\overline{F}^p_{i}-\xi^h_{pk}\overline{F}^p_{i}
 +\psi_i\delta^p_k\overline{F}^h_{p}+\psi_k\delta^p_i\overline{F}^h_{p}+\xi^p_{ik}\overline{F}^h_{p}\\
 &=\overline{F}^h_{i\underset1{\overline{|}}k}-\psi_p\delta^h_k\overline{F}^p_{i}-\psi_k\overline{F}^h_{i}-\xi^h_{pk}\overline{F}^p_{i}
 +\psi_i\overline{F}^h_{k}+\psi_k\overline{F}^h_{i}+\xi^p_{ik}\overline{F}^h_{p}
 \\
 &=\overset{0}{\overbrace{\overline{F}^h_{i\underset1{\overline{|}}k}}}-\psi_p\delta^h_k\overline{F}^p_{i}
 +\psi_i\overline{F}^h_{k}-\xi^h_{pk}\overline{F}^p_{i}+\xi^p_{ik}\overline{F}^h_{p},
\endaligned\end{eqnarray}
where $|$, and $\overline{|}$ are respectively covariant derivatives
in $\mathbb{GR}_N$ and $\mathbb{G}\overline{\underset1{\mathbb{K}}}{}_N$.
\end{demo}

\begin{te} The generalized Riemannian space
$\mathbb{GR}_N$ admits a nontrivial geodesic mapping onto the generalized
K\"ahlerian space $\mathbb{G}\overline{\underset1{\mathbb{K}}}{}_N$ if and only if, in
the common coordinate system $x$ with respect to the mapping,  the
conditions
\begin{equation}\label{mik}
\begin{array}{crl}
a)\;& \overline{g}_ {{{ij}\underset2|k}}{}=&\overline{g}_
{{\underset\lor{ij}\underset2{\overline{|}}k}}{}
+2\psi_k\overline{g}_{{ij}}+\psi_i\overline{g}_{{jk}}+\psi_j\overline{g}_{{ik}}+\xi^\alpha_{ki}\overline{g}_{\alpha
j}
+\xi^\alpha_{kj}\overline{g}_{i\alpha};\\
 b)&\;\overline{F}^h_{i\underset2|k}=&\overline{F}^h_{k}\psi_i -\delta^h_k\overline{F}^\alpha_{i}\psi_\alpha
 -\xi^h_{k\alpha }\overline{F}^\alpha_{i}+\xi^\alpha_{ki}\overline{F}^h_{\alpha};
\end{array}%
\end{equation}
hold with respect to the second kind of covariant derivatives, where $\psi_i\ne0$ and the tensors $\overline{g}_{{ij}}$ and $\overline{F}^h_{i}$ satisfy the following conditions:
\begin{equation}\label{uslov1}det(\overline{g}_{\underline{ij}})\ne0,\quad
\overline{F}^\alpha_{i}\overline{g}_{\underline{\alpha
j}}+\overline{F}^\alpha_{j}\overline{g}_{\underline{\alpha i}}=0,
\quad
\overline{F}^h_{\alpha}\overline{F}^\alpha_{i}=-\delta^h_i.\end{equation}
Then $\overline{g}_{{ij}}$ and $\overline{F}^h_{i}$ are
the metric tensor and the structure of
$\mathbb{G}\overline{\underset1{\mathbb{K}}}{}_N$, respectively.
 \end{te}

\begin{demo} In $\mathbb{GR}_N$ for the second kind of covariant
derivatives, we have
\begin{eqnarray}\aligned
\overline{F}^h_{i\underset2|k}&=\overline{F}^h_{i,k}+(\overline{\Gamma}^h_{kp}-\psi_k\delta^h_p-\psi_p\delta^h_k-\xi^h_{kp})
 \overline{F}^p_{i}-(\overline{\Gamma}^p_{ki}-\psi_k\delta^p_i-\psi_i\delta^p_k-\xi^p_{ki})\overline{F}^h_{p}\\
 &=\overline{F}^h_{i,k}+\overline{\Gamma}^h_{kp}\overline{F}^p_{i}-\psi_k\delta^h_p\overline{F}^p_{i}
 -\psi_p\delta^h_k\overline{F}^p_{i}-\xi^h_{kp}\overline{F}^p_{i}
 -\overline{\Gamma}^p_{ki}\overline{F}^h_{p}+\psi_k\delta^p_i\overline{F}^h_{p}+\psi_i\delta^p_k\overline{F}^h_{p}+\xi^p_{ki}\overline{F}^h_{p}\\
 &=\overline{F}^h_{i\underset2{\overline{|}}k}-\psi_k\delta^h_p\overline{F}^p_{i}-\psi_p\delta^h_k\overline{F}^p_{i}-\xi^h_{kp}\overline{F}^p_{i}
 +\psi_k\delta^p_i\overline{F}^h_{p}+\psi_i\delta^p_k\overline{F}^h_{p}+\xi^p_{ki}\overline{F}^h_{p}\\
 &=\overline{F}^h_{i\underset2{\overline{|}}k}-\psi_k\overline{F}^h_{i}-\psi_p\delta^h_k\overline{F}^p_{i}-\xi^h_{kp}\overline{F}^p_{i}
 +\psi_k\overline{F}^h_{i}+\psi_i\overline{F}^h_{k}+\xi^p_{ki}\overline{F}^h_{p}
 \\
 &= \psi_i\overline{F}^h_{k}-\psi_p\delta^h_k\overline{F}^p_{i}
  -\xi^h_{kp}\overline{F}^p_{i}+\xi^p_{ki}\overline{F}^h_{p}.
\endaligned\end{eqnarray}
\end{demo}

\smallskip
In the similar way, we can prove corresponding theorems for the third and the fourth kind of covariant derivative:

\begin{te} The generalized Riemannian space
$\mathbb{GR}_N$ admits a nontrivial geodesic mapping onto the generalized
K\"ahlerian space $\mathbb{G}\overline{\underset1{\mathbb{K}}}{}_N$ if and only if, in
the common coordinate system $x$ with respect to the mapping,  the
conditions
\begin{equation}\label{mik}
\begin{array}{crl}
a)\;& \overline{g}_ {{{ij}\underset3|k}}{}=&\overline{g}_
{{\underset\lor{ij}\underset3{\overline{|}}k}}{}
+2\psi_k\overline{g}_{{ij}}+\psi_i\overline{g}_{{jk}}+\psi_j\overline{g}_{{ik}}+\xi^\alpha_{ik}\overline{g}_{\alpha
j}
+\xi^\alpha_{kj}\overline{g}_{i\alpha};\\
 b)&\;\overline{F}^h_{i\underset3|k}=&\overline{F}^h_{i\underset3{\overline{|}}k} -\psi_p\delta^h_k\overline{F}^p_{i}
  +\psi_i\overline{F}^h_{k}-\xi^h_{pk}\overline{F}^p_{i}+\xi^p_{ki}\overline{F}^h_{p},
\end{array}%
\end{equation}
hold with respect to the third kind of covariant derivatives, where $\psi_i\ne0$ and the tensors $\overline{g}_{{ij}}$ and $\overline{F}^h_{i}$ satisfy the following conditions:
\begin{equation}\label{uslov1}det(\overline{g}_{\underline{ij}})\ne0,\quad
\overline{F}^\alpha_{i}\overline{g}_{\underline{\alpha
j}}+\overline{F}^\alpha_{j}\overline{g}_{\underline{\alpha i}}=0,
\quad
\overline{F}^h_{\alpha}\overline{F}^\alpha_{i}=-\delta^h_i.\end{equation}
Then $\overline{g}_{{ij}}$ and $\overline{F}^h_{i}$ are
the metric tensor and the structure of
$\mathbb{G}\overline{\underset1{\mathbb{K}}}{}_N$, respectively.
 \end{te}


\begin{te} The generalized Riemannian space
$\mathbb{GR}_N$ admits a nontrivial geodesic mapping onto the generalized
K\"ahlerian space $\mathbb{G}\overline{\underset1{\mathbb{K}}}{}_N$ if and only if, in
the common coordinate system $x$ with respect to the mapping,  the
conditions
\begin{equation}\label{mik}
\begin{array}{crl}
a)\;& \overline{g}_ {{{ij}\underset4|k}}{}=&\overline{g}_
{{\underset\lor{ij}\underset4{\overline{|}}k}}{}
+2\psi_k\overline{g}_{{ij}}+\psi_i\overline{g}_{{jk}}+\psi_j\overline{g}_{{ik}}+\xi^\alpha_{ki}\overline{g}_{\alpha
j}
+\xi^\alpha_{jk}\overline{g}_{i\alpha};\\
 b)&\;\overline{F}^h_{i\underset4|k}=&\overline{F}^h_{i\underset4{\overline{|}}k} -\psi_p\delta^h_k\overline{F}^p_{i}
  +\psi_i\overline{F}^h_{k}-\xi^h_{kp}\overline{F}^p_{i}+\xi^p_{ik}\overline{F}^h_{p},
\end{array}%
\end{equation}
hold with respect to the fourth kind of covariant derivatives, where $\psi_i\ne0$ and the tensors $\overline{g}_{{ij}}$ and $\overline{F}^h_{i}$ satisfy the following conditions:
\begin{equation}\label{uslov1}det(\overline{g}_{\underline{ij}})\ne0,\quad
\overline{F}^\alpha_{i}\overline{g}_{\underline{\alpha
j}}+\overline{F}^\alpha_{j}\overline{g}_{\underline{\alpha i}}=0,
\quad
\overline{F}^h_{\alpha}\overline{F}^\alpha_{i}=-\delta^h_i.\end{equation}
Then $\overline{g}_{{ij}}$ and $\overline{F}^h_{i}$ are
the metric tensor and the structure of
$\mathbb{G}\overline{\underset1{\mathbb{K}}}{}_N$, respectively.
 \end{te}

\subsection{Equitorsion geodesic mapping}

Equitorsion mappings play an important role in the theories of geodesic, conformal and holomorphically projective transformations between two
spaces of non-symmetric affine connection.

\begin{de} {\bf \cite{bib-65}} A  mapping $f:\mathbb{GR}_N\to \mathbb{G}\overline {\underset1{\mathbb{K}}}{}_N$ is
an {\bf equitorsion geodesic mapping} if the  torsion tensor of the
spaces $\mathbb{GR}_N$ and $\mathbb{G}\overline {\underset1{\mathbb{K}}}{}_N$ are equal. Then
from (\ref{0.2}), (\ref{0.3}) and (\ref{def})
\begin{equation}
\overline{\Gamma}{}^h_{\underset\lor {ij}}-\Gamma^h_{\underset\lor
{ij}}=\KX hij =0,\label{0.8}
\end{equation}
where $\underset\lor{ij}$ denotes an antisymmetrization with respect
to $i,j$.
\end{de}
In the case of these mappings, the previous Theorems 2.3.-2.6.
become:

\begin{te} The generalized Riemannian space
$\mathbb{GR}_N$ admits a nontrivial equitorsion geodesic mapping onto
the generalized K\"ahlerian space $\mathbb{G}\overline{\underset1{\mathbb{K}}}{}_N$ if and
only if, in the common coordinate system $x$ with respect to the
mapping,  the conditions
\begin{equation}\label{ }
\begin{array}{crl}
a)\;& \overline{g}_ {{{\underline{ij}}\underset1|k}}{}=&
2\psi_k\overline{g}_{{\underline{ij}}}+\psi_i\overline{g}_{{\underline{jk}}}+\psi_j\overline{g}_{{\underline{ik}}};\\
b)&\;\overline{F}^h_{i\underset1|k}=&\overline{F}^h_{k}\psi_i
-\delta^h_k\overline{F}^p_{i}\psi_p;
\end{array}%
\end{equation}
hold with respect to the first kind of covariant derivatives, where $\psi_i\ne0$ and the tensors $\overline{g}_{{ij}}$ and $\overline{F}^h_{i}$ satisfy the following conditions:
\begin{equation}\label{uslov1}det(\overline{g}_{\underline{ij}})\ne0,\quad
\overline{F}^\alpha_{i}\overline{g}_{\underline{\alpha
j}}+\overline{F}^\alpha_{j}\overline{g}_{\underline{\alpha i}}=0,
\quad
\overline{F}^h_{\alpha}\overline{F}^\alpha_{i}=-\delta^h_i.\end{equation}
 \end{te}

\medskip
\begin{te} The generalized Riemannian space
$\mathbb{GR}_N$ admits a nontrivial equitorsion geodesic mapping onto the
generalized K\"ahlerian space $\mathbb{G}\overline{\underset1{\mathbb{K}}}{}_N$ if and
only if, in the common coordinate system $x$ with respect to the
mapping,  the conditions
\begin{equation}\label{ }
\begin{array}{crl}
a)\;& \overline{g}_ {{{\underline{ij}}\underset2|k}}{}=&
2\psi_k\overline{g}_{{\underline{ij}}}+\psi_i\overline{g}_{{\underline{jk}}}+\psi_j\overline{g}_{{\underline{ik}}};\\
b)&\;\overline{F}^h_{i\underset2|k}=&\overline{F}^h_{k}\psi_i
-\delta^h_k\overline{F}^p_{i}\psi_p;
\end{array}%
\end{equation}
hold with respect to the second kind of covariant derivatives,
where $\psi_i\ne0$ and the tensors $\overline{g}_{{ij}}$ and $\overline{F}^h_{i}$ satisfy the following conditions:
\begin{equation}\label{uslov1}det(\overline{g}_{\underline{ij}})\ne0,\quad
\overline{F}^\alpha_{i}\overline{g}_{\underline{\alpha
j}}+\overline{F}^\alpha_{j}\overline{g}_{\underline{\alpha i}}=0,
\quad
\overline{F}^h_{\alpha}\overline{F}^\alpha_{i}=-\delta^h_i.\end{equation}
 \end{te}

\medskip
\begin{te} The generalized Riemannian space
$\mathbb{GR}_N$ admits a nontrivial equitorsion geodesic mapping onto
the generalized K\"ahlerian space $\mathbb{G}\overline{\underset1{\mathbb{K}}}{}_N$ if and
only if, in the common coordinate system $x$ with respect to the
mapping,  the conditions
\begin{equation}\label{ }
\begin{array}{crl}
a)\;& \overline{g}_ {{{\underline{ij}}\underset3|k}}{}=&2\psi_k\overline{g}_{{\underline{ij}}}
+\psi_i\overline{g}_{{\underline{jk}}}+\psi_j\overline{g}_{{\underline{ik}}};\\
 b)\;&\overline{F}^h_{i\underset3|k}=&\overline{F}^h_{i\underset3{\overline{|}}k} -\psi_p\delta^h_k\overline{F}^p_{i}
  +\psi_i\overline{F}^h_{k},
\end{array}%
\end{equation}
hold with respect to the third kind of covariant derivatives,
where $\psi_i\ne0$ and the tensors $\overline{g}_{{ij}}$ and $\overline{F}^h_{i}$ satisfy the following conditions:
\begin{equation}\label{uslov1}det(\overline{g}_{\underline{ij}})\ne0,\quad
\overline{F}^\alpha_{i}\overline{g}_{\underline{\alpha
j}}+\overline{F}^\alpha_{j}\overline{g}_{\underline{\alpha i}}=0,
\quad
\overline{F}^h_{\alpha}\overline{F}^\alpha_{i}=-\delta^h_i.\end{equation}
 \end{te}

\medskip
\begin{te} The generalized Riemannian space
$\mathbb{GR}_N$ admits a nontrivial equitorsion geodesic mapping onto the
generalized K\"ahlerian space $\mathbb{G}\overline{\underset1{\mathbb{K}}}{}_N$ if and
only if, in the common coordinate system $x$ with respect to the
mapping,  the conditions
\begin{equation}\label{ }
\begin{array}{crl}
a)\;& \overline{g}_ {{{\underline{ij}}\underset4|k}}{}=&
2\psi_k\overline{g}_{{\underline{ij}}}+\psi_i\overline{g}_{{jk}}+\psi_j\overline{g}_{{\underline{ik}}};\\
 b)\;&\overline{F}^h_{i\underset4|k}=&\overline{F}^h_{i\underset4{\overline{|}}k} -\psi_p\delta^h_k\overline{F}^p_{i}
  +\psi_i\overline{F}^h_{k},
\end{array}%
\end{equation}
hold with respect to the fourth kind of covariant derivatives,
where $\psi_i\ne0$ and the tensors $\overline{g}_{{ij}}$ and $\overline{F}^h_{i}$ satisfy the following conditions:
\begin{equation}\label{uslov1}det(\overline{g}_{\underline{ij}})\ne0,\quad
\overline{F}^\alpha_{i}\overline{g}_{\underline{\alpha
j}}+\overline{F}^\alpha_{j}\overline{g}_{\underline{\alpha i}}=0,
\quad
\overline{F}^h_{\alpha}\overline{F}^\alpha_{i}=-\delta^h_i.\end{equation}
 \end{te}

\section{Conclusion}

We continued the general idea by introducing the notion of generalized K\"ahlerian spaces of the
first kind $\mathbb{G}{\underset1 {\mathbb{K}}}{}_N$.

The main results of the paper:

1. New explicit formulas of geodesic mappings onto
$\mathbb{G}{\underset1 {\mathbb{K}}}{}_N$ are given in Section
2.

2. New explicit formulas of equitorsion geodesic mappings onto
$\mathbb{G}{\underset1 {\mathbb{K}}}{}_N$   are given in Section
2.1.

In a similar way we can consider generalized K\"ahlerian spaces of the second, the third and the fourth kind.

\end{document}